\documentclass[11pt]{article}
\usepackage{bigints}
\usepackage{color, amsmath,amssymb, amsfonts, amstext,amsthm, latexsym}

\usepackage{amssymb, epsfig, amssymb, latexsym}
\usepackage{amsmath}
\usepackage{graphicx}
\usepackage{longtable}
\usepackage{color}

  \newtheorem{Main Results}[theorem]{MainResults}

\newcommand{\f}{\mathbf f}
\newcommand{\g}{\mathbf g}
\newcommand{\y}{\mathbf y}

\begin{document}

\title{Modeling nonlinear random vibration: Implication of the energy conservation law
\footnote{This work was supported by NSFC grants 10971125 and 11028102, the NSF grant DMS-1025422, and the fundamental research funds for the Central Universities HUST 2010ZD037. } }

 \author{Xu Sun\\
 School of Mathematics and Statistics\\Huazhong University of Science and Technology\\ Wuhan 430074, Hubei, China \\
   E-mail: xsun15@gmail.com\\ \\
Jinqiao Duan\\
Institute for Pure and Applied Mathematics\\University of California, Los Angeles, CA 90095, USA\\ \&\\
Department of Applied Mathematics\\ Illinois Institute of Technology,
   Chicago, IL 60616, USA \\
   E-mail: duan@iit.edu\\ \\
    Xiaofan Li \\
Department of Applied Mathematics\\ Illinois Institute of Technology,   Chicago, IL 60616, USA \\   E-mail:    lix@iit.edu\\ \\
   }

\date{May 25, 2012}
\maketitle

\newpage
\begin{abstract}
Nonlinear random vibration under excitations of both Gaussian and Poisson white noises is considered.  The model is based on stochastic differential equations, and the corresponding stochastic integrals are defined in such a way that the energy conservation law is satisfied. It is shown that Stratonovich integral and Di Paola-Falsone integral should be used for excitations of Gaussian and Poisson white noises, respectively, in order for the model to satisfy the underlining physical laws (e.g., energy conservation). Numerical examples are presented to illustrate the theoretical results.

\bigskip


\end{abstract}

\section {Introduction}

Differential equations are extensively used in modeling dynamical systems in science and engineering. When dynamical systems are  under
random influences, stochastic differential equations (SDEs) may be more appropriate for modeling. The solutions of SDEs are   interpreted   in terms of stochastic integrals \cite{okse2003, klebaner2005}.

Dynamical  systems subject to Gaussian  white noise are often modeled by SDEs with Brownian motion, and the solutions are in terms of  the   Ito integral \cite{CWSTo2000, okse2003, klebaner2005}. Although the Ito integral is self-consistent mathematically, it is not the only type of stochastic integrals that can be constructed to interpret an SDE. Other
stochastic integrals, such as the Stratonovich integral \cite{CWSTo2000, okse2003, klebaner2005}, have also been used to interpret  an SDE as a stochastic integral equation. There is no right or wrong choice when choosing either Ito or
Stratonovich integrals in interpreting   SDEs mathematically, since the two integrals are equivalent and can be converted into each other, provided that the integrand satisfies certain smoothness conditions \cite{CWSTo2000, okse2003, klebaner2005}. However, these stochastic integrals have different definitions, and one may be more directly related to a practical situation
than the other. While Ito integral is a reasonable choice in many applications including finance and biology \cite{okse2003}, Stratonovich
integral is believed to be more appropriate in physical and engineering applications \cite{CWSTo2000}. Stratonovich integral has an extra term comparing with the corresponding Ito integral: the so-called   correction term \cite{CWSTo2000, WongZakai1965}. Some authors \cite{Ibrahim1985, YongLin1987} attribute this correction term to the conversion from physical white noise to ideal white noise. This explanation is not   necessarily convincing \cite{CWSTo2000, CWSTo1988}.

Dynamical systems driven by non-Gaussian white noise, especially Poisson white noise,  have    attracted a lot of attention  recently. Correction terms for converting Ito SDEs to Stratonovich SDEs with Poisson white noise   are presented in \cite{DiPaolaFalsone1993, DiPaolaFalsone1993b}. Although these correction terms have been accepted widely, there are some confusions \cite{Grigoriu1998, Hu1994}.

%

 In this paper, we  consider nonlinear random vibration under excitations of either Gaussian or Poisson white noises, modeled by appropriate stochastic differential equations. The main objective of this paper is to explain the correction terms in both Gaussian and Poisson white noise cases, from a physical perspective. We will show that the correction terms are natural consequences of   fundamental physical laws satisfied by the vibration system. Note that conventional spectral analysis methods \cite{LinCai2004},  which have found extensive applications in random vibration analysis, are not applicable in this case due to the nonlinearity of the system.

To this end, we consider a vibration system as a mass-spring-damping oscillator with random excitation
\begin{align}\label{NO}
m\overset{..}{ x}(t) + k  x(t) =g(x(t), \dot x(t)) + f(x(t), \dot x(t))\dot L(t),
\end{align}
where $m$ represents the mass, $k$ is the stiffness coefficient of the spring, $x(t)$ is the displacement depending on time $t$, and $\dot x(t)$ is the velocity. $ g(x(t), \dot x(t))$ and $ f(x(t), \dot x(t))\dot L(t)$ represent the generalized force terms, which may originate from external or parametric excitations. $\dot L(t)$ is a noise term defined as the formal derivative of some  stochastic process
\begin{align}\label{impulsivemodel}
L(t)= b  B(t) +c C(t),
\end{align}
where $b$ and $c$ are constants, $B(t)$ is a Gaussian process, and $C(t)$ is some compound Poisson process, which is expressed as
\begin{align}\label{section1_tmpp1}
C(t)=\sum_{i=1}^{N(t)} R_i U(t-t_i).
\end{align}
In Eq. (\ref{section1_tmpp1}), $N(t)$ is a Poisson process with intensity parameter $\lambda$, $U(t-t_i)$ is a unit step function (a Heaviside function) at $t_i$, $R_i$ is a random variable representing the $i$-th impulse. It follows from (\ref{impulsivemodel}) that
 \begin{align}\label{impulsivemodel_1}
\dot L(t)= b  \dot B(t) +c \dot C(t),
\end{align}
 where $\dot B(t)$ is the Gaussian white noise, and $\dot C(t)$ is the Poisson white noise expressed as
\begin{align}\label{dpn}
\dot C(t) = \sum_{i=1}^{N(t)} R_i\delta(t-t_i).
\end{align}
Note that (\ref{impulsivemodel_1}) expresses a general noise model including the Gaussian white noise ($b\ne 0, c=0$), the Poisson white noise ($b=0, c\ne 0$), and the combined Gaussian and Poisson white noise ($b\ne 0$ and $c \ne 0$).

The second-order equation (\ref{NO}) can be rewritten as a system of SDEs
\begin{align}\label{solutionofito_1}
d{\begin{pmatrix} x(t) \\\dot x(t) \end{pmatrix}} =\begin{pmatrix}0 & 1\\ -\dfrac{k}{m} & 0\end{pmatrix} \begin{pmatrix} x(t)\\\dot x(t) \end{pmatrix} \,{\rm d}t + \dfrac{1}{m} \begin{pmatrix} 0\\g(x(t), \dot x(t))\end{pmatrix} \,{\rm d}t+ \dfrac{1}{m} \begin{pmatrix} 0\\f(x(t), \dot x(t))\end{pmatrix} \,{\rm d}L(t).
\end{align}
Since $L(t)$ is non-differentiable almost everywhere, (\ref{solutionofito_1}) cannot be interpreted in the framework of classical calculus. Thus the solution of (\ref{solutionofito_1}) is interpreted with a stochastic integral,
\begin{align}\label{section1_tmp1}
 {\begin{pmatrix} x(t) \\\dot x(t) \end{pmatrix}}& =  {\begin{pmatrix} x(0) \\\dot x(0) \end{pmatrix}} + \int_0^t  \begin{pmatrix}0 & 1\\ -\dfrac{k}{m} & 0\end{pmatrix} \begin{pmatrix} x(s)\\\dot x(s) \end{pmatrix} \,{\rm d}s
    + \dfrac{1}{m} \int_0^t  \begin{pmatrix} 0\\g(x(s), \dot x(s))\end{pmatrix} \,{\rm d}s  \nonumber\\
    &\quad + \dfrac{1}{m} \int_0^t \begin{pmatrix} 0\\f(x(s), \dot x(s))\end{pmatrix} \,{\rm d}L(s).
\end{align}
Defining $\y (t) =\begin{pmatrix}x(t)\\\dot x(t) \end{pmatrix}$, $A=\begin{pmatrix} 0 & 1\\ -\dfrac{k}{m} & 0\end{pmatrix}$ and using the variation of parameters formula, the solution to Eq.~(\ref{NO}) can also be rewritten as \cite{okse2003},
\begin{align}\label{solutionofito}
  \y (t) = e^{A t} \y_0 +  \dfrac{1}{m} \int_0^t e^{A(t-s)}  \begin{pmatrix} 0\\g(x(s), \dot x(s))\end{pmatrix} \,{\rm   d} s + \dfrac{1}{m} \int_0^t e^{A(t-s)}\begin{pmatrix} 0\\f(x(s), \dot x(s))\end{pmatrix}\,{\rm d} L(s),
\end{align}
where $\y_0=\begin{pmatrix}x(0)\\\dot x(0) \end{pmatrix}\equiv \begin{pmatrix}x_0\\\dot x_0 \end{pmatrix}$ is the initial condition.   It can be shown that (\ref{section1_tmp1}) and (\ref{solutionofito}) are equivalent \cite{okse2003}.

It is straightforward to verify that
\begin{align}\label{MatrixA}
e^{A t} =I+\frac{At}{1!}+\frac{A^2 t^2}{2!}+\cdots =  \begin{pmatrix} \cos (\omega t) & \frac{\sin(\omega t)}{\omega}\\ -\omega  \sin(\omega t) & \cos (\omega t)\end{pmatrix},
\end{align}
where  $\omega = \sqrt{\dfrac{k}{m}}$. Substituting (\ref{MatrixA}) into (\ref{solutionofito}), we get
\begin{align}\label{solutionstobe}
\begin{cases}
x (t)& = \cos(\omega t) x_0 +\dfrac{\sin (\omega t)}{\omega} \dot x_0
  +  \bigints_0^t  \dfrac{\sin \omega(t-s)}{m\omega} g(x(s),\dot x(s))\, {\rm d}s \\
&\quad + \bigints_0^t \dfrac{\sin\omega(t-s) }{m\omega} f(x(s),\dot x(s))\, {\rm d}L(s),\\
\dot x (t)& = -\omega \sin (\omega t) x_0 +\cos (\omega t) \dot x_0 + \bigints_0^t \dfrac{\cos \omega (t-s)}{m} g(x(s),\dot x(s)) \,{\rm d}s\\
&\quad + \bigints_0^t \dfrac{\cos \omega (t-s)}{m} f(x(s),\dot x(s)) \,{\rm d}L(s).
\end{cases}
\end{align}
\normalsize

Note that the stochastic integrals in Eqs.~(\ref{section1_tmp1}) and (\ref{solutionstobe}) are yet to be defined. As stated earlier,  the stochastic integrals which can be used to interpreted SDEs may not be unique.  The question is that which stochastic integral will lead to the solution that is consistent with the physics of the system. One possible answer is to compare solutions to the SDEs with the corresponding experimental results. However, this method may be impractical in many cases due to the high cost of performing the experiments, as one needs highly accurate data from sufficiently large number of samples in order to resolve the subtle difference in the theory. In this paper, to construct a SDE model that is physically relevant to the real system, we propose to apply a stochastic integral such that the fundamental physical law (e.g., energy conservation) is satisfied.

This paper is organized as follows. In Sec.~\ref{sec.2}, starting from the energy conservation law, we define the stochastic integral that is suitable for the SDE model of the nonlinear random oscillators. The relationship between the proposed models and the existing models is discussed in Sec.~\ref{sec.rel}. Numerical methods with an illustrative example are presented in Sec.~\ref{sec.nm}.

\section{Stochastic integrals  for nonlinear oscillators under noise excitation} \label{sec.2}

Since there are multiple forms of the stochastic integrals that can be constructed from the SDE, we define the stochastic integral such that the fundamental physical laws are satisfied.  As for the nonlinear oscillators described by the SDE (\ref{NO}), we expect the energy-work conservation be satisfied
\begin{align}\label{energyprincipal_0}
&\left[\frac{1}{2}m\dot x^2(t) + \frac{1}{2}k x^2(t)\right]- \left[\frac{1}{2}m\dot x^2(0) + \frac{1}{2}k x^2(0) \right]\nonumber\\
&\quad = \int_0^t  \left[ g(x(s), \dot x(s)) + f(x(s), \dot x(s))\dot L(s)\right] \,{\rm d}x(s),
\end{align}
where $\frac{1}{2}m\dot x^2(t) + \frac{1}{2}k x^2(t)$ represents the total mechanical energy of the system at time $t$, and the integrand in the right hand side is the forcing term of (\ref{NO}). Equation (\ref{energyprincipal_0}) expresses that the change in the total mechanical energy is equal to the work done by the external forces.

 Writing (\ref{energyprincipal_0}) in the form of stochastic integral, we have
\begin{align}\label{energyprincipal_1}
&\left[\frac{1}{2}m\dot x^2(t) + \frac{1}{2}k x^2(t)\right]\nonumber- \left[\frac{1}{2}m\dot x^2(0) + \frac{1}{2}k x^2(0) \right]\nonumber\\
&\quad =   \int_0^t  g(x(s), \dot x(s))\dot x(s) \,{\rm d}s +   \int_0^t f(x(s), \dot x(s)) \dot L (s) \dot x(s) \,{\rm d}s\nonumber\\
&\quad =   \int_0^t  g(x(s), \dot x(s))\dot x(s) \,{\rm d}s +   \int_0^t f(x(s), \dot x(s))   \dot x(s) \,{\rm d}L(s).
\end{align}

As stated before, the stochastic integral with respect to $L(t)$ should be defined such that the solution of (\ref{NO}) satisfies the energy conservation law (\ref{energyprincipal_1}).
It follows from (\ref{impulsivemodel}) that a stochastic integral with respect to $L(t)$ can be decomposed into two terms:  stochastic integral with respect to $B(t)$ and stochastic integral with respect to $C(t)$.
We define the two terms in the next two subsections.

\subsection{For Gaussian white noises}\label{sec.gwn}

Assume $b=1$ and $c=0$, then it follows from (\ref{impulsivemodel}) that the stochastic process $L(t)$  reduces to  a Brownian motion, and (\ref{solutionstobe}) and (\ref{energyprincipal_1}) become

\begin{align}\label{solutions_ttmmpp}
\begin{cases}
x (t)& = \cos(\omega t) x_0 +\dfrac{\sin (\omega t)}{\omega} \dot x_0 + \bigints_0^t \dfrac{\sin \omega(t-s)}{m\omega} g(x(s),\dot x(s)) \,{\rm d}s\\
 &\quad +  \bigints_0^t \dfrac{\sin\omega(t-s) }{m\omega} f(x(s),\dot x(s))\,{\rm d}B(s),\\
\dot x (t)& = -\omega \sin (\omega t) x_0 +\cos (\omega t) \dot x_0 + \bigints_0^t \dfrac{\cos \omega (t-s) }{m} g(x(s),\dot x(s))\,{\rm d}s \\
&\quad +  \bigints_0^t \dfrac{\cos \omega (t-s)}{m}  f(x(s),\dot x(s))\,{\rm d}B(s).
\end{cases}
\end{align}
\normalsize
and
\begin{align}\label{energy_tmp1}
&\left[\frac{1}{2}m\dot x^2(t) + \frac{1}{2}k x^2(t)\right]- \left[\frac{1}{2}m\dot x^2(0) + \frac{1}{2}k x^2(0) \right]\nonumber\\
&=     \int_0^t  g(x(s), \dot x(s))  \dot x(s) \,{\rm d}s
   +   \int_0^t  f(x(s), \dot x(s)) \dot x(s) \,{\rm d}B(s),
\end{align}
respectively.

There are two types of stochastic integral extensively used for SDEs driven by Brownian motions: Ito integral and Stratonovich integral. Throughout this paper, we use '$\star$' to denote Ito calculus, and '$\circ$' for Stratonovich calculus. In the sense of Ito, (\ref{solutions_ttmmpp}) and (\ref{energy_tmp1}) can be written as

\begin{align}\label{solutions_ito}
\begin{cases}
x (t) &= \cos(\omega t) x_0 +\dfrac{\sin (\omega t)}{\omega} \dot x_0
 +  \bigints_0^t \dfrac{\sin \omega(t-s)}{m\omega} g(x(s),\dot x(s))\,{\rm d}s\\
 &\quad +  \bigints_0^t \dfrac{\sin\omega(t-s) }{m\omega} f(x(s),\dot x(s))\star {\rm d}B(s),\\
\dot x (t)& = -\omega \sin (\omega t) x_0 +\cos (\omega t) \dot x_0 + \bigints_0^t \dfrac{\cos \omega (t-s)}{m} g(x(s),\dot x(s))\,{\rm d}s\\
&\quad +  \bigints_0^t \dfrac{\cos \omega (t-s)}{m} f(x(s),\dot x(s))\star {\rm d}B(s).
\end{cases}
\end{align}
\normalsize
and
\begin{align}\label{energyprincipal_ito}
&\left[\frac{1}{2}m\dot x^2(t) + \frac{1}{2}k x^2(t)\right]- \left[\frac{1}{2}m\dot x^2(0) + \frac{1}{2}k x^2(0) \right]\nonumber\\
&=     \int_0^t  g(x(s), \dot x(s))  \dot x(s) \,{\rm d}s
   +   \int_0^t  f(x(s), \dot x(s)) \dot x(s)\star \,{\rm d}B(s),
\end{align}
respectively.
In the sense of Stratonovich, (\ref{solutions_ttmmpp}) and (\ref{energy_tmp1}) can be written as
\begin{align}\label{solutions_stratonovich}
\begin{cases}
x (t)& = \cos(\omega t) x_0 +\dfrac{\sin (\omega t)}{\omega} \dot x_0
   + \bigints_0^t \dfrac{\sin \omega(t-s)}{m\omega} g(x(s),\dot x(s))\,{\rm d}s \\
&\quad + \bigints_0^t \dfrac{\sin\omega(t-s) }{m\omega} f(x(s),\dot x(s))\circ {\rm d}B(s),\\
\dot x (t)& = -\omega \sin (\omega t) x_0 +\cos (\omega t) \dot x_0
   + \bigints_0^t \dfrac{\cos \omega (t-s)}{m} g(x(s),\dot x(s))\,{\rm d}s\\
&+ \bigints_0^t \dfrac{\cos \omega (t-s)}{m} f(x(s),\dot x(s))\circ {\rm d}B(s).
\end{cases}
\end{align}
\normalsize
and
\begin{align}\label{energy_stratonovich}
&\left[\frac{1}{2}m\dot x^2(t) + \frac{1}{2}k x^2(t)\right]- \left[\frac{1}{2}m\dot x^2(0) + \frac{1}{2}k x^2(0) \right]\nonumber\\
&=   \int_0^t  g(x(s), \dot x(s))  \dot x(s) \,{\rm d}s
 +    \int_0^t  f(x(s), \dot x(s)) \dot x(s)\circ \,{\rm d}B(s),
\end{align}
respectively.
Provided that the function $f$ is sufficient smooth,
the solutions in Stratonovich integrals,
(\ref{solutions_stratonovich}) and (\ref{energy_stratonovich}),
can be converted into the following forms with Ito integrals \cite{okse2003}
\small
\begin{align}\label{solutions_stratonovich_2}
\begin{cases}
x(t)& =\cos(\omega t) x_0 +\dfrac{\sin (\omega t)}{\omega} \dot x_0
  + \bigints_0^t \dfrac{\sin \omega(t-s)}{m\omega} \left[ g(x(s),\dot x(s)) + \dfrac{1}{2m} f(x(s), \dot x(s))f_{\dot x} (x(s), \dot x(s))\right] \,{\rm d}s \\
& + \bigints_0^t \dfrac{\sin \omega(t-s)}{m\omega} f(x(s), \dot x(s))\star {\rm d}B(s),\\
\dot x(t)& =-\omega \sin(\omega t) x_0 +\cos(\omega t) \dot x_0
  +  \bigints_0^t \dfrac{\cos \omega (t-s)}{m} \left[g(x(s),\dot x(s))
  + \dfrac{1}{2m} f(x(s), \dot x(s)) f_{\dot x} (x(s), \dot x(s))\right]   \,{\rm d}s\\
&+ \bigints_0^t \dfrac{\cos \omega (t-s)}{m} f(x(s), \dot x(s)) \star {\rm d}B(s),
\end{cases}
\end{align}
\normalsize
and
\begin{align}\label{energy_stratonovich_2}
&\left[\frac{1}{2}m\dot x^2(t) + \frac{1}{2}k x^2(t)\right]- \left[\frac{1}{2}m\dot x^2(0) + \frac{1}{2}k x^2(0) \right]\nonumber\\
&=    \int_0^t  \left[\left(g(x(s),\dot x(s)) + \frac{1}{2m}f(x(s), \dot x(s))f_{\dot x} (x(s), \dot x(s))   \right)\dot x(s)+\frac{1}{2m}f^2(x(s), \dot x(s))  \right]   \,{\rm d}s\nonumber\\
&\quad  +     \int_0^t  f(x(s), \dot x(s)) \dot x(s)\star {\rm d}B(s).
\end{align}

As shown in the Appendix, the solution (\ref{solutions_stratonovich_2}) satisfies the energy-work relation~\eqref{energy_stratonovich_2}, suggesting that when the randomness is modeled in sense of Stratonovich, the energy-work conservation law is satisfied. On the other hand, in a similar procedure as in the Appendix, it can be shown that the energy-work law (\ref{energyprincipal_ito}) contradicts with the solution (\ref{solutions_ito}). Therefore, Stratonovich integral instead of Ito integral should be used so that this nonlinear random oscillator model satisfies the energy conservation law.

This implies that when  Gaussian noise  is present in this nonlinear vibration system, the SDE model should be interpreted in the sense of   Stratonovich stochastic integral, but not in the sense of Ito stochastic integral.

\subsection{For Poisson white noises}\label{sec.pwn}
When $b=0$ and $c=1$, the stochastic process as expressed in (\ref{impulsivemodel}) reduces to a compound Poisson process. Note that  the jump size of $C(s)$ at time $s$ can be expressed as $\Delta C(s)=C(s)- C(s-)$, where $C(s-)$ is the left limit of $C(s)$ at $s$. Suppose $C(s)$ jumps at times $t_i$ ($i=1, 2, \cdots$), then the solution (\ref{section1_tmp1}) can be written as
\begin{align}
\begin{cases}\label{section22_tmp2}
 x(t)=x(0)+ \bigints_0^t \dot x(s) \,{\rm d}s,\\
   \dot x(t) =\dot x(0) -\dfrac{k}{m} \bigints_0^t x(s) \,{\rm d}s + \dfrac{1}{m} \bigints_0^t  g(x(s), \dot x(s)) \,{\rm d}s + \dfrac{1}{m}  \sum\limits_{i=1}^{N(t)}\bigints_{t_i-}^{t_i} f(x(s), \dot x(s))  \,{\rm d}C(s),
 \end{cases}
\end{align}
where $N(t)$, as shown in (\ref{section1_tmpp1}), represents the number of jumps upto time $t$.

In the following, we shall derive the stochastic integral with respect to jumps  such that the energy conservation law is satisfied.
First, let's examine the changes in the system at $i$-th jump occured
at time $t_i$ ($ 1\le i\le N(t)$).
From (\ref{section22_tmp2}), the displacement $x$ is continuous
while the velocity $\dot x$ undergoes an jump given by
\begin{align}\label{section22_tt2}
\begin{cases}
 x(t_i)=x(t_i-),\\
   \dot x(t_i) =\dot x(t_i-) + \dfrac{1}{m} \bigints_{t_i-}^{t_i} f(x(s), \dot x(s))  \,{\rm d}C(s).
 \end{cases}
\end{align}
The change in the total energy (\ref{energyprincipal_1}) due to the $i$-th
jump is that in the kinetic energy given by
\begin{align}\label{section22_tmp3}
 \frac{1}{2}m\dot x^2 (t_i ) - \frac{1}{2}m\dot x^2 (t_i-)
=   \int_{t_i-}^{ t_i } f(x(s), \dot x(s)) \dot x(s) \,{\rm d}C(s),
\end{align}
due to the continuity of the displacement $x$ across an jump.

If the integrals with respect to jumps are defined in sense of Ito, then (\ref{section22_tt2}) and (\ref{section22_tmp3}) becomes
\begin{align}\label{1-May152012}
   \dot x(t_i) =\dot x(t_i-) + \dfrac{1}{m} f(x(t_i-), \dot x(t_i-))  \Delta C(t_i).
 \end{align}
 and
\begin{align}\label{2-May152012}
 \frac{1}{2}m \dot x^2 (t_i ) -  \frac{1}{2}m \dot x^2 (t_i-)
=  f(x(t_i-), \dot x(t_i-)) \dot x(t_i-) \Delta C(t_i),
\end{align}
respectively. Since $\Delta C(t_i) \ne 0$, it is obvious that (\ref{1-May152012}) contradict with (\ref{2-May152012}), which indicates that the energy conservation law cannot be satisfied when the integrals with respect to jumps are interpreted in sense of Ito. In the following, we shall show that the integrals should be interpreted as some kind of Riemann integral on the imaginary path along the jump to satisfy the energy conservation law.

Let $\underline {\dot x} (t_i, r)$ be the value of $\dot x(s)$ at time $t_i$ if $C(s)$ jumped from $C(t_i-)$ to $r$.  Then $\underline {\dot x} (t_i, C(t_i-)) = \dot x(t_i-)$ and $\underline {\dot x} (t_i, C(t_i))=\dot x(t_i)$. With the integrals being interpreted as the Riemann integral on the imaginary path along  the jump, the energy-work law (\ref{section22_tmp3}) can be written as
\begin{align}\label{section22_tt1}
  \frac{1}{2}m \underline {\dot x}^2 (t_i, C(t_i) ) - \frac{1}{2}m \underline {\dot x}^2 (t_i, C(t_i-))  =   \int_{C(t_i-)}^{ C(t_i )} f(x(t_i), \underline {\dot x}(t_i, r)) \underline {\dot x}(t_i, r) \, {\rm d}r.
\end{align}
and the solution (\ref{section22_tt2}) becomes
\begin{align}\label{3-May152012}
   \underline {\dot x} (t_i, C(t_i) ) -  \underline {\dot x}  (t_i, C(t_i-))  = \dfrac{1}{m} \int_{C(t_i-)}^{ C(t_i )} f(x(t_i), \underline {\dot x}(t_i, r)) \, {\rm d}r.
\end{align}
Since the jump size can be any value, it follows from (\ref{section22_tt1}) and (\ref{3-May152012})that for any $ \lambda \in R$, it is true that
\begin{align}\label{section22_tmp4}
  \frac{1}{2}m\underline {\dot x}^2 (t_i, \lambda ) - \frac{1}{2} m \underline {\dot x}^2 (t_i, C(t_i-))  =   \int_{C(t_i-)}^{ \lambda} f(x(t_i), \underline {\dot x}(t_i, r)) \underline {\dot x}(t_i, r) \,{\rm d}r,
\end{align}
and
\begin{align}\label{4-May152012}
   \underline {\dot x} (t_i, \lambda ) -  \underline {\dot x}  (t_i, C(t_i-))  = \dfrac{1}{m} \int_{C(t_i-)}^{ \lambda} f(x(t_i), \underline {\dot x}(t_i, r)) \, {\rm d}r.
\end{align}
Taking derivatives of both sides of (\ref{section22_tmp4}) and (\ref{4-May152012}) with respect to $\lambda$, respectively, we get the identical ordinary
differential equation(ODE)
\begin{align}\label{section22_tmp7}
 \frac{\rm d}{{\rm d}\lambda}\underline {\dot x}(t_i, \lambda) = \dfrac{1}{m} f(x(t_i),\underline {\dot x}(t_i, \lambda)).
\end{align}
Therefore, the energy conservation law is satisfied.

Using the fact that $\underline {\dot x} (t_i, C(t_i-)) = \dot x(t_i-)$ and $\underline {\dot x} (t_i, C(t_i )) = \dot x(t_i )$ , it follows from (\ref{section22_tmp7}) that
\begin{align}\label{section22_tmp5}
  \dot x(t_i ) = \dot x(t_i-) +  \dfrac{1}{m}  Y_i(\Delta C(t_i)),
\end{align}
where $ Y_i(\Delta C(t_i))$ is determined by the initial or terminal value
problem of the ODE
\begin{align}\label{section22_tt3}
\begin{cases}
 \dfrac{{\rm d}}{{\rm d}\lambda}Y_i(\lambda) =f(x(t_i),Y_i(\lambda)+\dot x(t_i-)),\\
 Y_i(0)=0.
 \end{cases}
\end{align}
Note that in (\ref{section22_tt3}), $0 \le \lambda \le \Delta C(t_i)$ for $\Delta C(t_i)>0$ or $ \Delta C(t_i)\le \lambda\le 0$ for $\Delta C(t_i)<0$.
Comparing  the original solution expression (\ref{section22_tt2}) with the new formula (\ref{section22_tmp5}), it can be seen that the last term in (\ref{section22_tt2}) should be defined as $ Y_i(\Delta C(t_i))$, and hence (\ref{section22_tmp2})  should be interpreted as
\begin{align}
\begin{cases}
x(t)=x(0) + \int_0^t   \dot x(s) \,{\rm d}s,  \\
  \dot x(t ) = \dot x(0) +  \dfrac{1}{m} \bigints_0^t g(x(s), \dot x(s)) \,{\rm d}s + \dfrac{1}{m}  \sum\limits_{i=1}^{N(t)}  Y_i(\Delta C(t_i)),
  \end{cases}
\end{align}
where $Y_i( \Delta C(t_i))$ is the solution to the ODE (\ref{section22_tt3}).

As it will become clear in Sec.~\ref{sec.rel}, this implies that when pure jump noise, such as Poisson white noise, is present in this nonlinear vibration system, the SDE model should be interpreted in the sense of Di Paola-Falsone stochastic integral \cite{DiPaolaFalsone1993, DiPaolaFalsone1993b}.

\subsection{For combined Gaussian and Poisson white noises}
When both $b\ne 0$ and $c\ne 0$, the excitation is a combined Gaussian and Poisson white noise.
Combining the results in the subsections \ref{sec.gwn} and \ref{sec.pwn},
we find that, in order to satisfy the energy-work conservtion law,
one has to interpret the stochastic integrals with respect to Brownian motions as Stratonovich integrals, and the integrals with respect to jumps as DiPaola-Falsone integrals.
Therefore, the solution to \eqref{NO} is given by the expression
\eqref{section1_tmp1}, where the stochastic integral is defined as
\begin{align}\label{section23_tmp2b}
\int_0^t f(x(s), \dot x(s)) \,{\rm d}L(s) = b\int_0^t f(x(s), \dot x(s)) \circ \,{\rm d}B_s +  \sum\limits_{i=1}^{N(t)}  Y_i(  \Delta L(t_i)),
\end{align}
where $\displaystyle{\Delta L(t_i) = c [C(t_i)-C(t_i-)]=c \Delta  C(t_i)}$.
Recall that, in (\ref{section23_tmp2b}),
'$\circ$' denote integrals in the Stratonovich sense,
$N(t)$ represents the number of jumps up to time $t$,
and $Y_i( \Delta L(t_i))$ is the solution to the ODE (\ref{section22_tt3}),
where $\lambda$ takes value of $0 \le \lambda \le \Delta L(t_i)$ for $\Delta L(t_i)>0$ or $ \Delta L(t_i)\le \lambda\le 0$ for $\Delta L(t_i)<0$.

\section{Relationship with the existing models}\label{sec.rel}

In this section, we shall  show that the correction term $Y_i(\Delta L(t_i)) $, as given by the solution to the ODE (\ref{section22_tt3}),
is consistent with the one proposed in the work by Di Paola and Falsone
\cite{DiPaolaFalsone1993, DiPaolaFalsone1993b}.

When $f$ is Lipschitz continuous, it is easy to check that the solution of  (\ref{section22_tt3}) exists and is unique. If $f$ is assumed to be smooth, then $Y_i(\Delta L(t_i))$ is analytic with respect to $\Delta L(t_i)$.

Using Taylor expansion,
\begin{align}\label{section3_3}
Y_i(\Delta L(t_i)) =Y_i (0) + \frac{{ \rm  d}}{{\rm d}\lambda}  Y_i (\lambda)\big |_{\lambda=0}(\Delta L(t_i)) +\frac{1}{2!} \frac{{{ \rm d}^2}}{{\rm d}\lambda^2}  Y_i (\lambda) \big |_{\lambda=0}(\Delta L(t_i))^2+\cdots.
\end{align}
It follows from  (\ref{section22_tt3}) that for any $n\ge 1$,
\begin{align}\label{section3_3ttp}
\frac{{\rm d}^n }{{\rm d}\lambda^n}Y_i(\lambda) = \frac{{\rm d}}{{\rm d} Y_i(\lambda)}\left\{ \frac{{\rm d}^{n-1} }{{\rm d}\lambda^{n-1}} Y_i(\lambda)\right\}  f(x(t_i), Y_i(\lambda)+\dot x(t_i-)).
\end{align}
Substituting  (\ref{section3_3ttp}) into (\ref{section3_3}), and using the fact that $Y_i(0)=0$, we get
\begin{align}\label{section3_7}
Y_i(\Delta L(t_i)) =  \sum_{j=1}^{\infty}
\frac{ f^{(j)} (x(t_i), \dot x(t_i-))}{j!} \left(\Delta L(t_i)\right) ^j
\end{align}
where
\begin{align}\label{section3_5}
\begin{cases}
f^{(1)}(x(t), \dot x(t_i-)) =f(x(t_i),\dot x(t_i-)) \\
f^{(j)}(x(t), \dot x(t_i-)) =\left.\dfrac{\partial f^{(i-1)} (x(t_i),\dot x(t_i-)+\lambda)}{\partial \lambda}\right.\Big|_{\lambda=0}f(x(t_i),\dot x(t_i-)) &\text{ for } j\geq 2.
\end{cases}
\end{align}
Thus the correction term given by (\ref{section3_7}) is exactly the same as the one proposed in \cite{DiPaolaFalsone1993, DiPaolaFalsone1993b}.

We have shown that the correction term $Y_i(\Delta L(t_i))$ can be obtained
in two ways: solving the initial value problem
to the ODE (\ref{section22_tt3}) or computing the expansion (\ref{section3_7}).
Note that the former approach of solving (\ref{section22_tt3})
is applicable under much more general condition
than the latter one of evaluating the infinite series (\ref{section3_7}),
because the existence of the solution to the ODE only
requires $f(x,y)$ to be integrable but (\ref{section3_7}) demands $f(x,y)$ to be infinitely differentiable.

\section{Simulation examples} \label{sec.nm}

Solutions of the SDE for a nonlinear oscillator (\ref{NO}), defined  by (\ref{section1_tmp1}) and (\ref{section23_tmp2b}),  can hardly be obtained with analytical methods. In this section, the SDE (\ref{NO}) is numerically solved to verify the conclusion obtained in section 2. Consider the case with both Gaussian and Poisson white noises,
${\rm d}L = b\,{\rm d}B + c\,{\rm d}C$, with the compound Poisson
process given by \eqref{dpn}.

The numerical procedure of the SDE for a nonlinear oscillator (\ref{NO}) defined  by (\ref{section1_tmp1}) and (\ref{section23_tmp2b})  is as follows. On each time subinterval
$t_{i-1} < t < t_i$,
i.e. when no jumps occur, (\ref{NO}) becomes
 \begin{align}\label{section3_tmp1}
\begin{cases}
 {\rm d}x(t) =\dot x(t),  \,{\rm d}t\\
 {\rm d} \dot x(t)  =-\dfrac{k}{m} x (t) \,{\rm d}t + \dfrac{1}{m} g(x(t), \dot x(t))   \,{\rm d}t + \dfrac{b}{m} f(x(t), \dot x(t)) \circ \,{\rm d}B.
 \end{cases}
 \end{align}
The above equation can be converted into Ito SDE and then computed by conventional algorithms for Ito SDEs, such as Euler method, Milstein method, or other algorithm of high-order accuracy based on stochastic Taylor expansion \cite{KloedenPlaten1992}.
At the time $t_i$ when an jump occurs, it follows from (\ref{section1_tmp1})
and (\ref{section23_tmp2b}) that
  \begin{align}\label{section3_tmp2}
\begin{cases}
 x(t_i)= x(t_i-),\\
  \dot x (t_i)  =\dot x(t_i-) + \frac{1}{m} Y_i(\Delta L (t_i)),
 \end{cases}
 \end{align}
where $Y_i(\Delta L(t_i))$ is obtained by solving the deterministic
ODE (\ref{section22_tt3}) using Runge-Kutta or multistep methods.

Consider the following stochastic Duffing-van der Pol equation
\begin{align}\label{example_e1}
\overset{..}x+(1+  x^2 )\dot x + x + x^3 = \dot x \dot L(t),
\end{align}
with the initial condition $x(0)=2$ and $\dot x(0)=0$.
The SDE (\ref{example_e1}) can be written in form of (\ref{NO}) with $m=1$, $k=1$, $g(x, \dot x)= -(1+  x^2 )\dot x  -  x^3$ and $f(x, \dot x)=\dot x$.
In the simulation, we take $L(t)$ as in \eqref{impulsivemodel}
with $b=c=1$, i.e. $L(t) = B(t) + C(t)$,
where $C(t)$ is a pure jump process given by \eqref{section1_tmpp1}
with $N(t)$ being a Poisson process with intensity parameter as $\lambda=3.4$
and $R_i$ ($i=1, \cdots, N(t)$) being   random numbers of the standard normal distribution.

\bigskip
\noindent {\bf Case 1}

\begin{figure}[scale=1]
  \epsfig{file=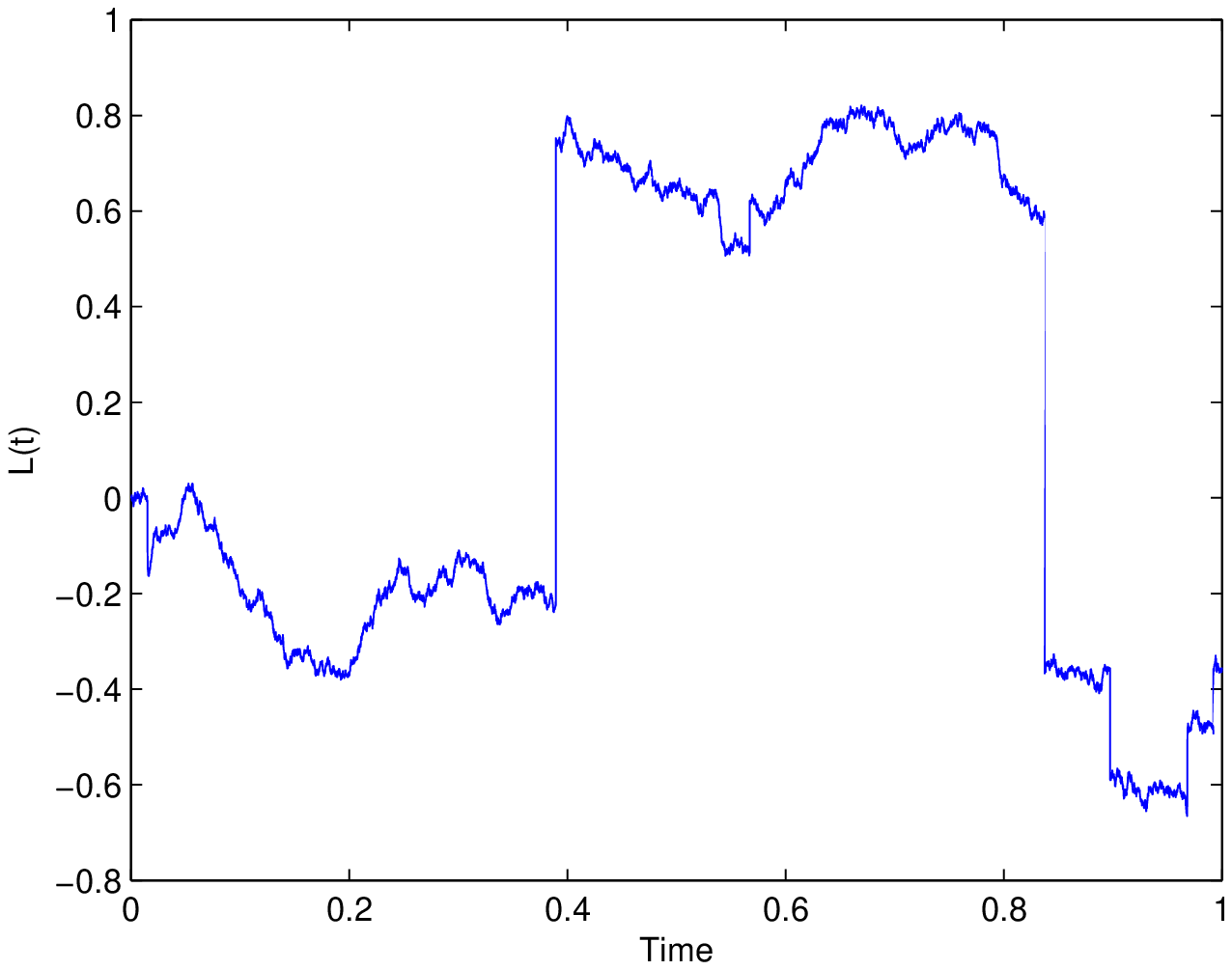,width=\linewidth}
  \caption{A sample path of the driving process $L(t)$ as a combination
of a Gaussian process and a compound Poisson process given in \eqref{section1_tmpp1}.}
  \label{fig0}
\end{figure}

\begin{figure}
  \epsfig{file=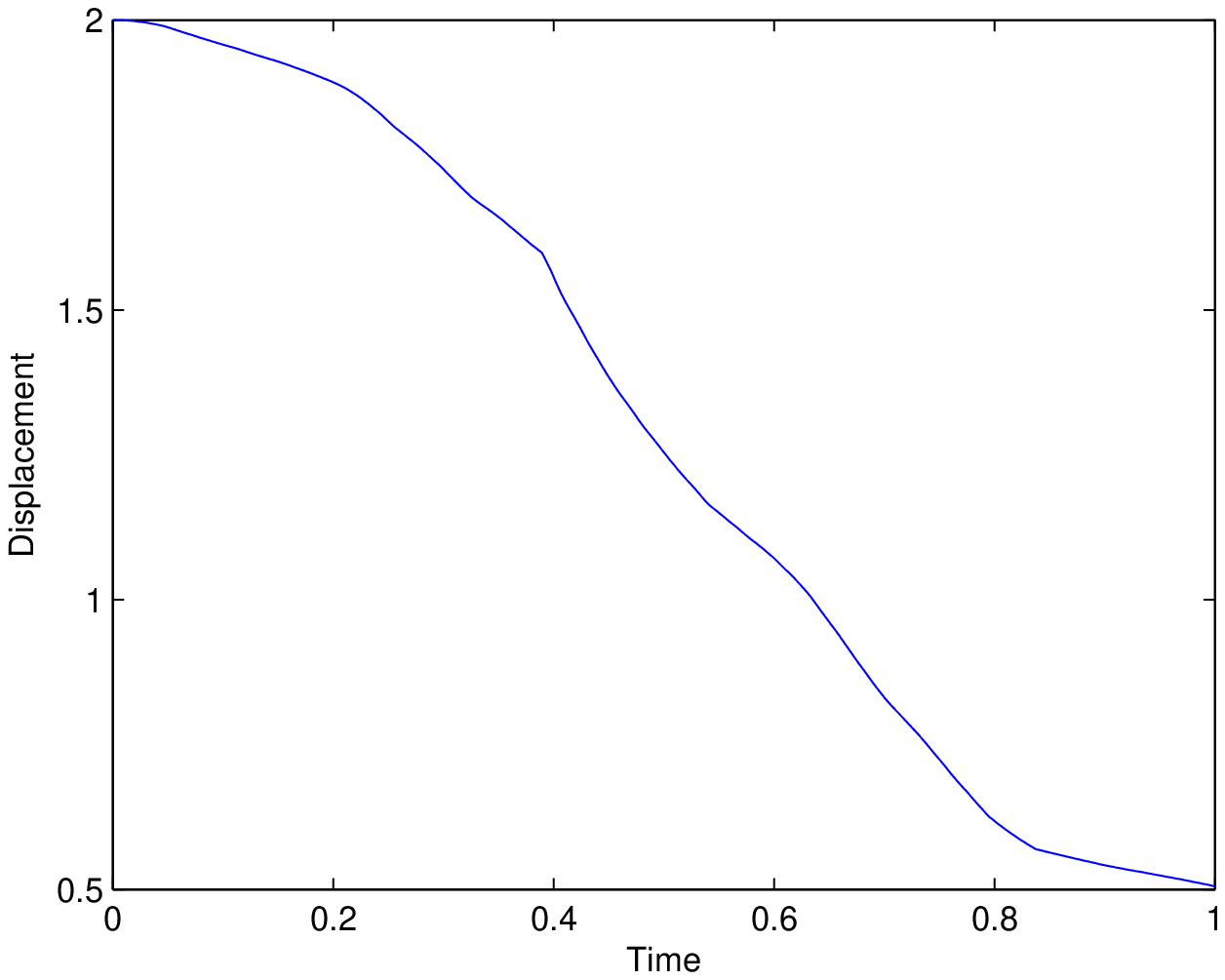,width=\linewidth}
  \caption{The evolution of the displacement $x(t)$ as the solution of Duffing-van der Pol equation \eqref{example_e1} defined by Stratonovich integral and Di Paola-Falsone integral, where the driving process $L(t)$ is given in Fig.~\ref{fig0}. }
  \label{fig1}
\end{figure}

\begin{figure}
  \epsfig{file=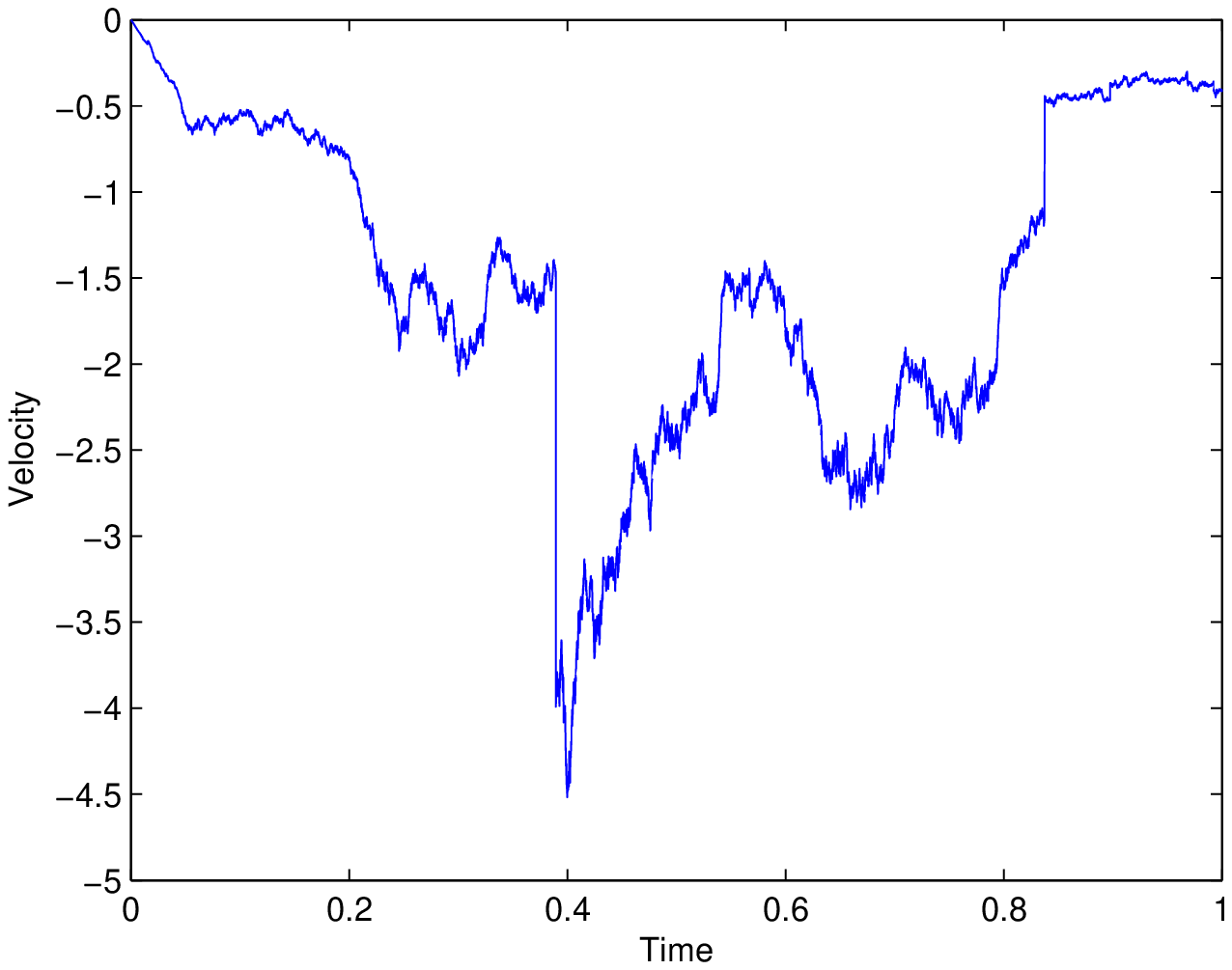,width=\linewidth}
  \caption{The evolution of the velocity $\dot x(t)$ as the solution of Duffing-van der Pol equation \eqref{example_e1} defined by Stratonovich integral and Di Paola-Falsone integral, where the driving process $L(t)$ is given in Fig.~\ref{fig0}. }
  \label{fig2}
\end{figure}

\begin{figure}
  \epsfig{file=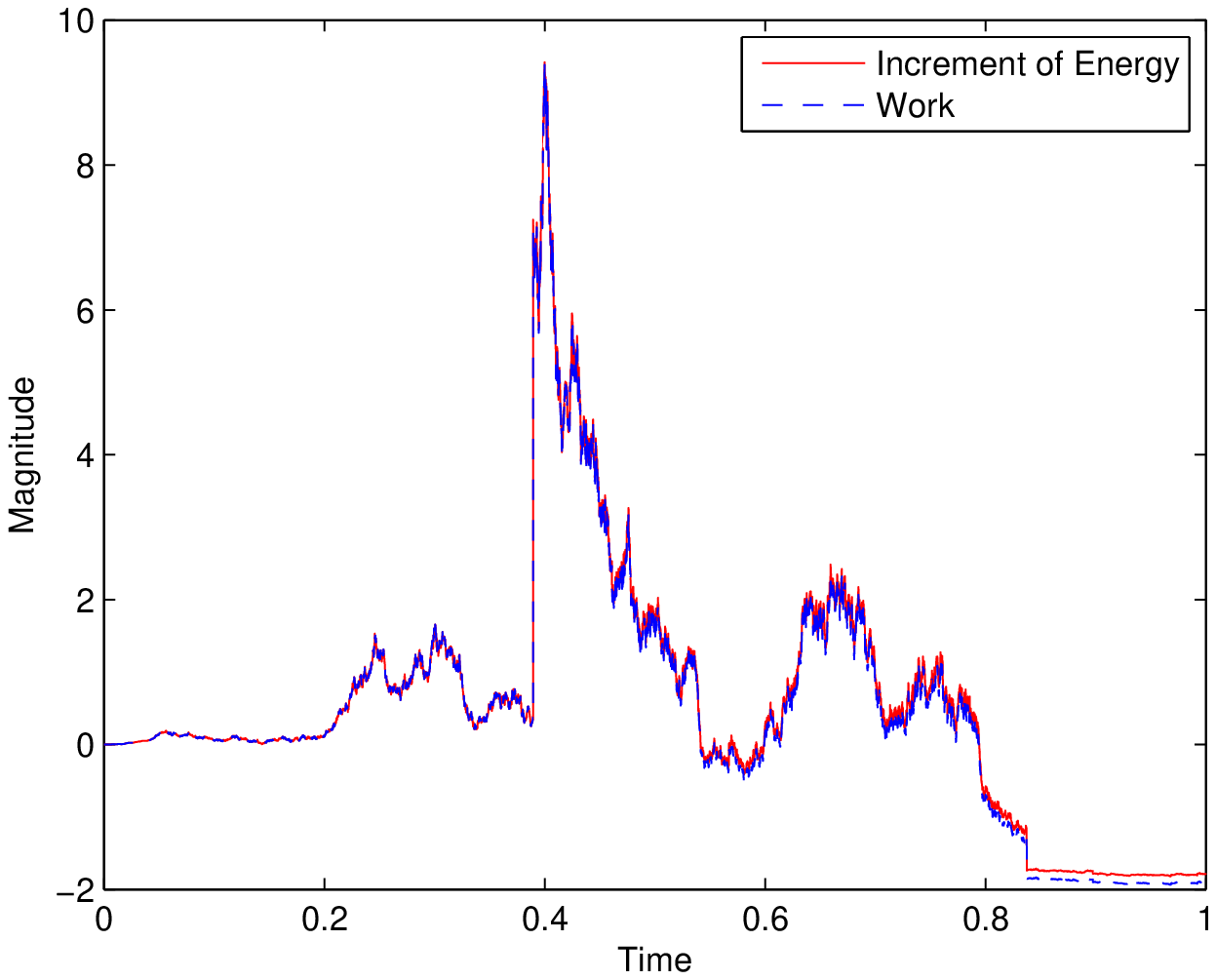,width=\linewidth}
  \caption{Comparison of the change in total energy, defined by (\ref{May222012-1}), and the work done by the force, defined by (\ref{May222012-2}).}
  \label{fig3}
\end{figure}

In this case, (\ref{example_e1}) or (\ref{NO}) is interpreted by (\ref{section1_tmp1}) and (\ref{section23_tmp2b}). In the simulation, to integrate (\ref{section1_tmp1}) and (\ref{section23_tmp2b}) numerically, we use Euler's method to advance (\ref{section3_tmp1}) when no jumps occur, while evaluate (\ref{section3_tmp2}) by solving (\ref{section22_tt3}) with Euler's method when jumps arrive. Note that to apply Euler's methods,  the stochastic integral in (\ref{section3_tmp1}) need to convert into Ito integral.
The step size of Euler's method for solving both (\ref{section3_tmp1}) and (\ref{section22_tt3}) is $\Delta t=0.0001$.
Figure~\ref{fig0} shows a sample path of the driven process $L(t)$, and Figs. \ref{fig1} and \ref{fig2} show the numerical solution of the
displacement $x(t)$ and the velocity $\dot x(t)$ respectively, corresponding
to the path shown in Fig.~\ref{fig0}.
Based on the numerical solutions shown in Figs.~\ref{fig1} and \ref{fig2},
 we compare  in Fig.~\ref{fig3} the energy increment $\Delta \text{E}(t)$,
 which is defined by
 \begin{align}\label{May222012-1}
 \Delta E =  \dfrac{1}{2} \left[x^2(t)+\dot x^2(t)-x^2(0)-\dot x^2(0)\right],
 \end{align}
 with the work done $\text{WK}(t)$
defined by $\displaystyle{\int_0^t  g(x, \dot x) \dot x(s)\,{\rm d}s}$
$\displaystyle{+ \int_0^t f(x, \dot x) \dot x(s) \,{\rm d}L(s)}$ where the stochastic
integral is taken in the sense similar to \eqref{section23_tmp2b}, i.e.,
\begin{align}\label{May222012-2}
WK(t)& = \int_0^t  g(x, \dot x) \dot x(s)\,{\rm d}s  + \int_0^t f(x, \dot x) \dot x(s) \,{\rm d}L(s)\nonumber\\
& = \int_0^t  g(x, \dot x) \dot x(s)\,{\rm d}s  + \int_0^t f(x, \dot x) \dot x(s)\circ \,{\rm d}B_s +   \sum\limits_{i=1}^{N(t)}  \bar Y_i(  \Delta L(t_i)),
\end{align}
where $\displaystyle{\Delta L(t_i) = c [C(t_i)-C(t_i-)]=c \Delta  C(t_i)}$,
'$\circ$' denote integrals in the Stratonovich sense,
$N(t)$ is the number of jumps upto time $t$,
and $\bar Y_i( \Delta L(t_i))$ is the solution to the following ODE,
\begin{align}\label{May222012-3}
\begin{cases}
 \dfrac{{\rm d}}{{\rm d}\lambda} \bar Y_i(\lambda) =f(x(t_i),\bar Y_i(\lambda)+\dot x(t_i-))(\bar Y_i(\lambda)+\dot x(t_i-)),\\
 \bar Y_i(0)=0.
 \end{cases}
\end{align}
with $\lambda$ taking value of $0 \le \lambda \le \Delta L(t_i)$ for $\Delta L(t_i)>0$ or $ \Delta L(t_i)\le \lambda\le 0$ for $\Delta L(t_i)<0$.
It can be seen clearly from Fig.~\ref{fig3} that the energy increment of the system agrees with the work very well,
indicating that the energy conservation law is satisfied.

\bigskip
\noindent {\bf Case 2}

\begin{figure}
  \epsfig{file=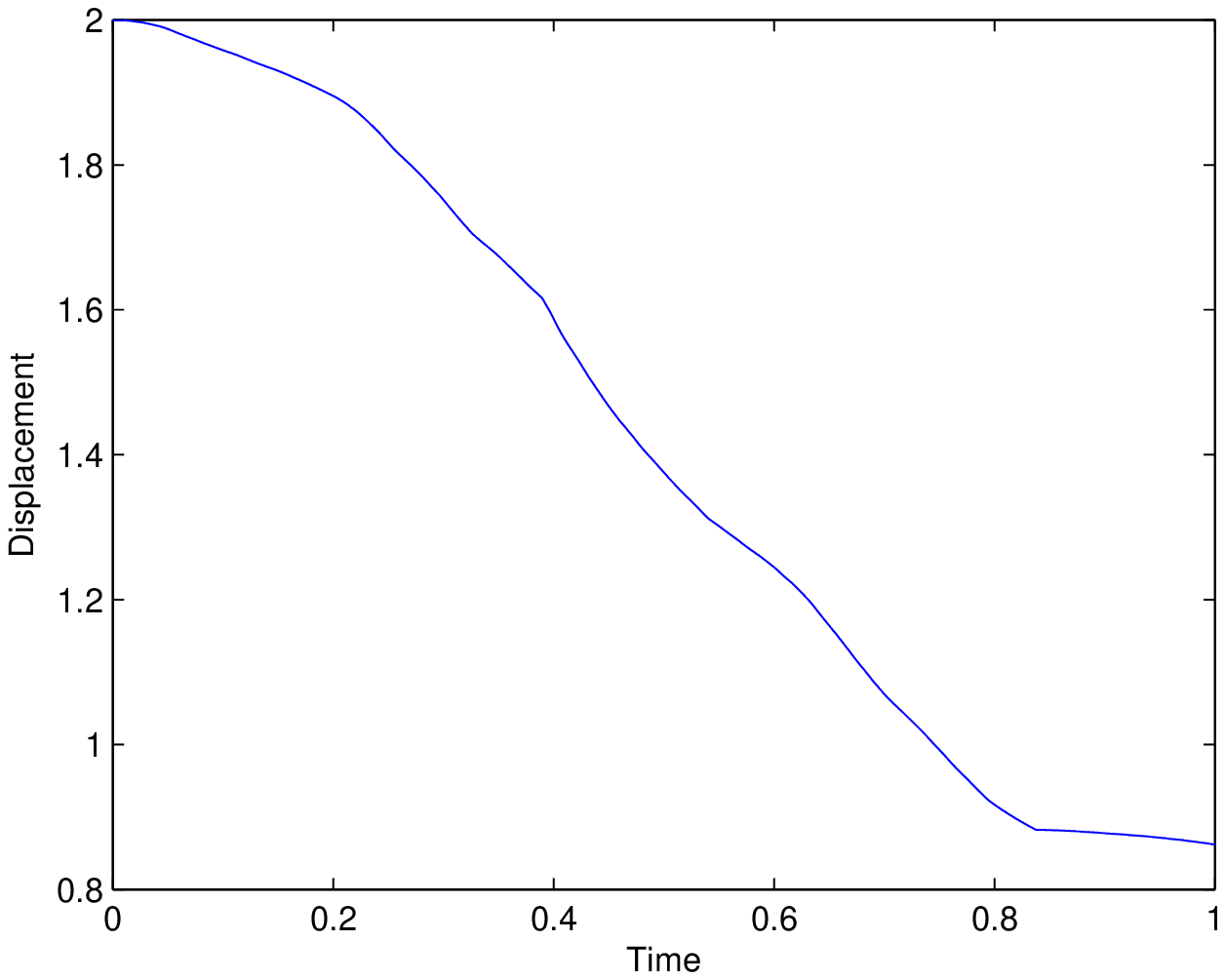,width=\linewidth}
  \caption{The evolution of the displacement $x(t)$ as the solution of Duffing-van der Pol equation \eqref{example_e1} defined by Ito integral, where the driving process $L(t)$ is given in Fig.~\ref{fig0}. }
  \label{fig5}
\end{figure}

\begin{figure}
  \epsfig{file=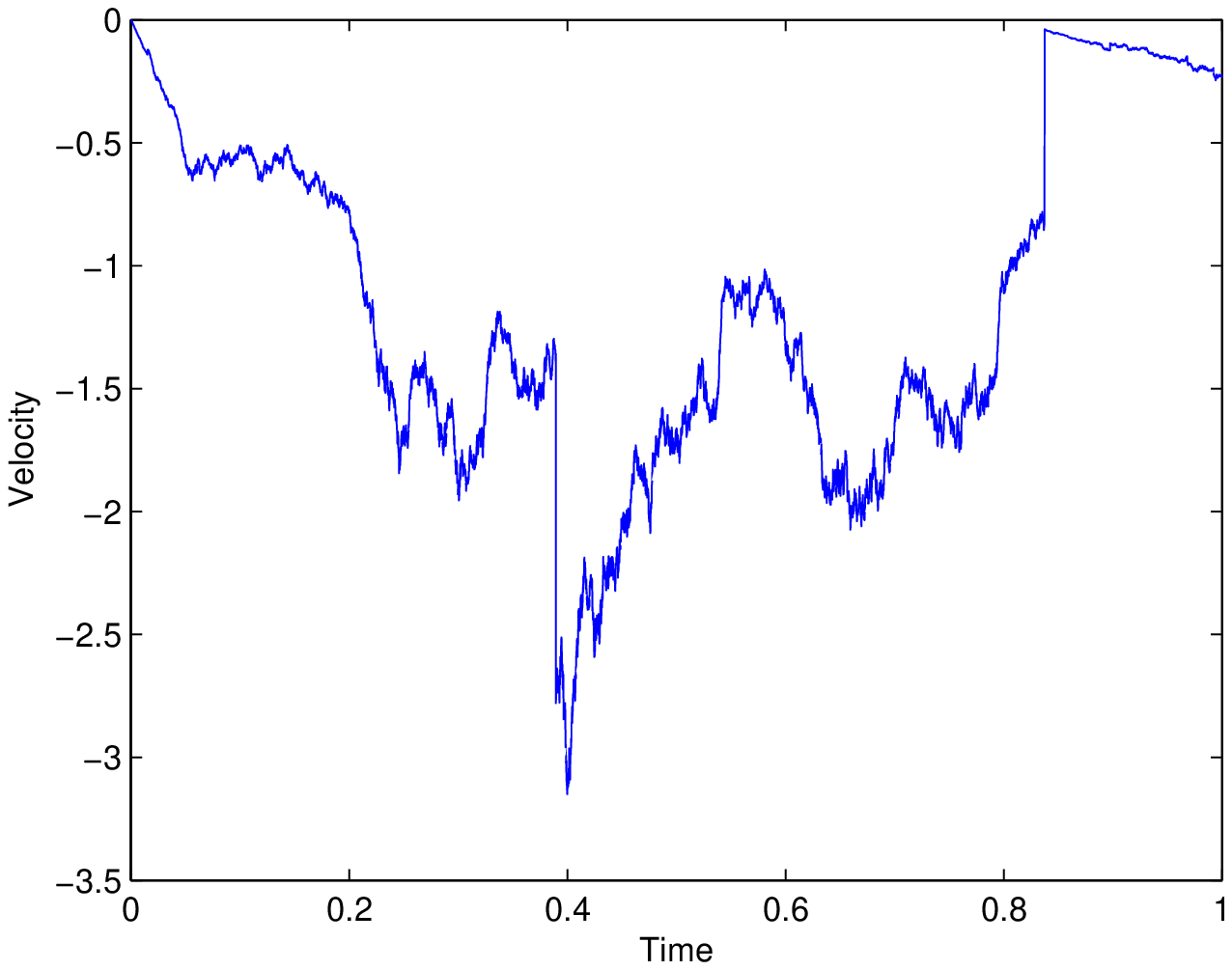,width=\linewidth}
  \caption{The evolution of the velocity $\dot x(t)$ as the solution of Duffing-van der Pol equation \eqref{example_e1} defined by Ito integral, where the driving process $L(t)$ is given in Fig.~\ref{fig0}. }
  \label{fig6}
\end{figure}

\begin{figure}
  \epsfig{file=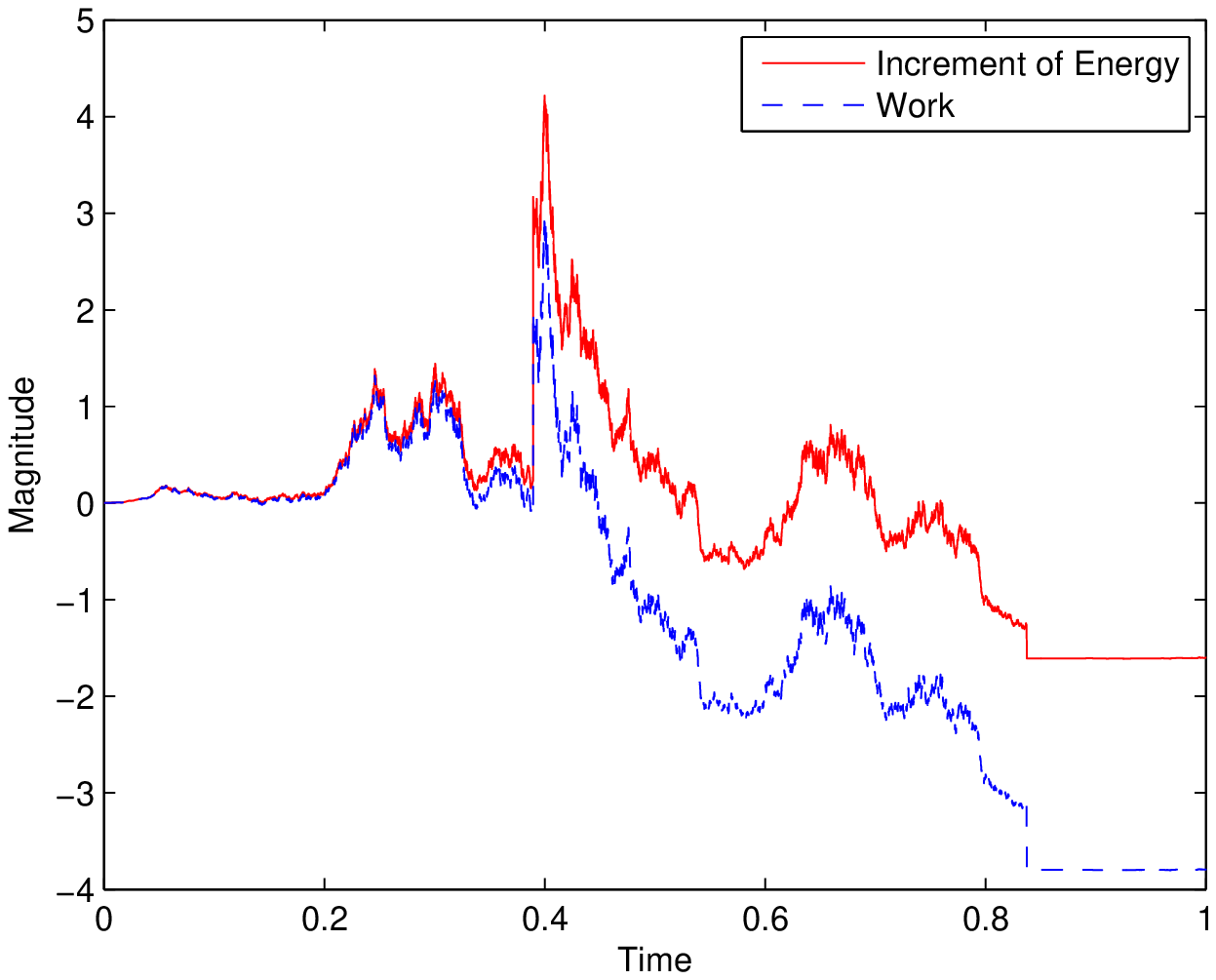,width=\linewidth}
  \caption{Comparison of the change in total energy, defined by (\ref{May222012-1}), and the work done by the force, defined by (\ref{May222012-6}). }
  \label{fig7}
\end{figure}

In this case, (\ref{example_e1}) (or (\ref{NO}) is  interpreted by using Ito stochastic integrals.
Now (\ref{section3_tmp1}) and (\ref{section3_tmp2}) become
  \begin{align}\label{May222012-4}
\begin{cases}
 {\rm d}x(t) =\dot x(t),  \,{\rm d}t\\
 {\rm d} \dot x(t)  =-\dfrac{k}{m} x (t) \,{\rm d}t + \dfrac{1}{m} g(x(t), \dot x(t))   \,{\rm d}t + \dfrac{b}{m} f(x(t), \dot x(t)) \star \,{\rm d}B,
 \end{cases}
 \end{align}
 and
  \begin{align}\label{May222012-5}
\begin{cases}
 x(t_i)= x(t_i-),\\
  \dot x (t_i)  =\dot x(t_i-) +\frac{1}{m} f(x(s-), \dot x(s-)) \Delta L(t_i),
 \end{cases}
 \end{align}
 respectively.
In the simulation, the driving process $L(t)$  and all the simulation parameters are taken the same as in case 1. Figures ~\ref{fig5} and ~\ref{fig6} present the corresponding numerical solutions of the displacement and velocity, respectively. Comparison of the energy increment, defined by (\ref{May222012-1}),  and the work done, now defined by
 \begin{align}\label{May222012-6}
WK(t) = \int_0^t  g(x, \dot x) \dot x(s)\,{\rm d}s   + b \int_0^t f(x, \dot x) \dot x(s)\star \,{\rm d}B_s +   \sum\limits_{i=1}^{N(t)}  f(x(s-), \dot x(s-))\dot x(s-)\Delta L(t_i) ,
\end{align}
 is presented in Fig.~\ref{fig7}. We can see clearly from Fig.~\ref{fig7} that there is significant difference between the energy increment and the work done by the force. Note that all the curves  in Fig.~\ref{fig7}  tend  to have a very small variance in the time span $0.84<t \le 1$. This is the consequence of the fact that  the velocity is very small for  $0.84<t \le 1$, as shown in Fig.~\ref{fig6}. Since $f(x(t),\dot x(t))=g(x(t),\dot x(t))=\dot x(t)$, it follows from (\ref{May222012-1}),  (\ref{May222012-4}) and (\ref{May222012-6}) that both the energy increment and the work done change slowly for  very small velocity $\dot x(t)$.

 \bigskip
 By comparing Fig.~\ref{fig7} with Fig.~\ref{fig3}, we can see that Stratonovich integral and Di Paola-Falsone integral should be used for excitations of Gaussian and Poisson white noises, respectively, in order for the model to satisfy the underlining physical laws.


%

\appendix
\section{Appendix: Proof of the energy-work law (\ref{energy_stratonovich_2})
from the solution (\ref{solutions_stratonovich_2})}

For cosmetic purpose, we introduce the following simplified notations: $g(s)=\dfrac{1}{m} g(x(s),\dot x(s))$, $f(s)=\dfrac{1}{m} f(x(s), \dot x(s))$, $f_{\dot x} (s)=\dfrac{1}{m} f_{\dot x(s)}(x(s), \dot x(s))$.
Moreover, in this Appendix, all the stochastic integrals
with respect to Brownian motions are in sense of Ito (we have dropped $\star$
notation).
Then the energy-work law (\ref{energy_stratonovich_2}) is equivalent to
\begin{align}\label{6-May152012}
& \frac{1}{2} \left[ \dot x^2(t) + \omega^2 x^2(t)\right]- \frac{1}{2} \left[\ \dot x^2(0) + \omega^2 x^2(0) \right]\nonumber \\
& =     \int_0^t  \left[\left(g(s) +  \dfrac{1}{2} f(s)f_{\dot x} (s)   \right)\dot x(s)+\frac{1}{2 }f^2(s)  \right]   \,{\rm d}s
  +     \int_0^t  f(s) \dot x(s)\, {\rm d}B(s).
\end{align}
Next, we show the solution given in (\ref{solutions_stratonovich_2})
satisfies the energy-work law (\ref{6-May152012}).

Denote the right and left hand sides of  (\ref{6-May152012}) as $RHS$ and $LHS$, respectively. Substitute (\ref{solutions_stratonovich_2}) into the left hand side of (\ref{6-May152012}),  we get
\small
\begin{align}\label{A1_1}
 &LHS =\frac{1}{2} \left( \int_0^t \sin \omega (t-s)\left[g(s) + \frac{1}{2}f(s) f_{\dot x} (s) \right] \,{\rm d}s \right)^2+\frac{1}{2} \left( \int_0^t \sin(\omega (t-s)) f(s) \,{\rm d}B(s) \right)^2\nonumber \\
&\quad +  \frac{1}{2} \left( \int_0^t \cos \omega (t-s)\left[g(s) + \frac{1}{2}f(s) f_{\dot x} (s) \right] \,{\rm d}s \right)^2+\frac{1}{2}\left( \int_0^t \cos(\omega (t-s)) f(s) \,{\rm d}B(s) \right)^2 \nonumber\\
&\quad+ ( \omega \cos(\omega t) x_0 + \sin(\omega t) \dot x_0)\int_0^t \sin \omega (t-s)\left[g(s) + \frac{1}{2}f(s) f_{\dot x} (s) \right] \,{\rm d}s\nonumber\\
& \quad + ( \omega \cos(\omega t) x_0 + \sin(\omega t) \dot x_0)\int_0^t \sin(\omega (t-s)) f(s)\,{\rm d}B(s)\nonumber\\
&\quad+  \left(\int_0^t \sin \omega (t-s) \left[g(s) + \frac{1}{2}f(s) f_{\dot x} (s) \right] \,{\rm d}s \right) \left( \int_0^t \sin(\omega(t-s)) f(s) \,{\rm d}B(s)\right) \nonumber\\
&\quad+ (-\omega \sin(\omega t) x_0 + \cos(\omega t) \dot x_0)\int_0^t \cos \omega (t-s)\left[g(s) + \frac{1}{2}f(s) f_{\dot x} (s) \right] \,{\rm d}s\nonumber\\
& \quad + (-\omega \sin(\omega t) x_0 + \cos(\omega t) \dot x_0)\int_0^t \cos(\omega (t-s)) f(s) \,{\rm d}B(s)\nonumber\\
&\quad+  \left( \int_0^t \cos(\omega (t-s) \left[g(s) + \frac{1}{2}f(s) f_{\dot x} (s) \right] \,{\rm d}s\right) \left(\int_0^t \cos(\omega(t-s)) f(s) \,{\rm d}B(s)\right).
\end{align}
\normalsize
Substituting (\ref{solutions_stratonovich_2}) into the right-hand side of (\ref{6-May152012}),  we get
\begin{align}\label{A1_2}
 &RHS = \int_0^t \left(g(s)+ f(s)f_{\dot x}(s) \right)\left(-\omega \sin(\omega s) x_0 + \cos(\omega s) \dot x_0\right) \,{\rm d}s\notag\\
 &\quad+ \int_0^t\int_0^s  \left[\cos(\omega(s-p)) \left(g(s)+  \frac{f(s)f_{\dot x}(s)}{2}\right) \left(g(p)+  \frac{f(p)f_{\dot x}(p)}{2}\right)\right]\,{\rm d}p\,{\rm d}s\notag\\
 &\quad+ \int_0^t\int_0^s  \left[\cos(\omega(s-p)) \left(g(s)+  \frac{f(s)f_{\dot x}(s)}{2}\right) f(p)\right]\,{\rm d}B(p) \,{\rm d}s  + \frac{1}{2}\int_0^t f^2(s) \,{\rm d}s\nonumber\\
  &\quad + \int_0^t f(s) \left(-\omega \sin(\omega s) x_0+ \cos(\omega s) \dot x_0 \right)\,{\rm d}B(s)\nonumber\\
 &\quad + \int_0^t \int_0^s \cos(\omega(s-p)) f(s)\left(g(p)+  \frac{f(p)f_{\dot x}(p)}{2}\right) \,{\rm d}p \,{\rm d}B(s)\notag\\
 &\quad + \int_0^t \int_0^s \cos(\omega (s-p)) f(s)f(p) \,{\rm d}B(p)\,{\rm d}B(s).
\end{align}
To prove $LHS$ in \eqref{A1_1} is equal to $RHS$ in \eqref{A1_2},
we claim the following facts
\begin{align}\label{A1_3}
&\int_0^t\int_0^s  \cos(\omega(s-p)) \left(g(s)+  \frac{f(s)f_{\dot x}(s)}{2}\right) \left(g(p)+  \frac{f(p)f_{\dot x}(p)}{2}\right) \,{\rm d}p\, {\rm d}s\notag\\
&=  \frac{1}{2} \left[ \int_0^t \sin(\omega (t-s)) \left(g(s)+  \frac{f(s)f_{\dot x}(s)}{2}\right) \,{\rm d}s \right]^2 \notag\\
 &\quad + \frac{1}{2}  \left[ \int_0^t \cos(\omega (t-s)) \left(g(s)+  \frac{f(s)f_{\dot x}(s)}{2}\right)\,{\rm d}s \right]^2,
\end{align}
\begin{align}\label{A1_4}
& \int_0^t \int_0^s \cos(\omega (s-p)) f(s)f(p) \,{\rm d}B(p) \,{\rm d}B(s) + \frac{1}{2} \int_0^t f^2(s) \,{\rm d}s\nonumber\\
&=\frac{1}{2} \left( \int_0^t \sin(\omega (t-s)f(s) \,{\rm d}B(s) \right)^2 +\frac{1}{2} \left( \int_0^t \cos(\omega (t-s)f(s) \,{\rm d}B(s) \right)^2,
\end{align}
\begin{align}\label{A1_5}
&\int_0^t\int_0^s  \left[\cos(\omega(s-p)\left(g(s)+  \frac{f(s)f_{\dot x}(s)}{2}\right) \right] f(p) \,{\rm d}B(p) \,{\rm d}s \nonumber\\
&\quad + \int_0^t\int_0^s  \left[\cos(\omega(s-p)\left(g(p)+  \frac{f(p)f_{\dot x}(p)}{2}\right) \right] f(s) \,{\rm d}p\,{\rm d}B(s)\notag\\
&=    \int_0^t \cos(\omega (t-s)) \left(g(s)+  \frac{f(s)f_{\dot x}(s)}{2}\right)\,{\rm d}s \int_0^t \cos(\omega(t-s)) f(s)  \,{\rm d}B(s) \notag\\
 &\quad +   \int_0^t \sin(\omega (t-s)) \left(g(s)+  \frac{f(s)f_{\dot x}(s)}{2}\right)\,{\rm d}s \int_0^t \sin(\omega(t-s))f(s) \,{\rm d}B(s)\,
\end{align}
\begin{align}\label{A1_6}
&\int_0^t \left(g(s)+  \frac{f(s)f_{\dot x}(s)}{2}\right) \left(-\omega \sin(\omega s) x_0 + \cos(\omega s) \dot x_0\right) \,{\rm d}s\nonumber\\
&=  (-\omega \sin(\omega t) x_0 + \cos(\omega t) \dot x_0)\int_0^t \cos(\omega (t-s))\left(g(s)+  \frac{f(s)f_{\dot x}(s)}{2}\right) \,{\rm d}s\nonumber\\
&\quad + (\omega\cos(\omega t) x_0 + \sin(\omega t) \dot x_0)\int_0^t \sin(\omega (t-s))\left(g(s)+  \frac{f(s)f_{\dot x}(s)}{2}\right) \,{\rm d}s,
\end{align}
and
\begin{align}\label{A1_7}
&\int_0^t f(s) \left(-\omega \sin(\omega s) x_0+ \cos(\omega s) \dot x_0 \right)\,{\rm d}B(s)\nonumber\\
&=  (-\omega \sin(\omega t) x_0 + \cos(\omega t) \dot x_0)\int_0^t \cos(\omega (t-s)) f(s) \,{\rm d}B(s)\nonumber\\
&\quad +  (\omega \cos(\omega t) x_0 + \sin(\omega t) \dot x_0)\int_0^t \sin(\omega (t-s)) f(s) \,{\rm d}B(s).
\end{align}
One can easily see that (\ref{A1_6}) and (\ref{A1_7}) are true by using
the trignometric identities
\begin{align}
\cos(\omega s)=\cos(\omega t)\cos(\omega (t-s)) +\sin(\omega t)\sin(\omega (t-s))\notag,
 \end{align}
 and
\begin{align}
\sin(\omega s)=\sin(\omega t)\cos(\omega (t-s)) -\cos(\omega t)\sin(\omega (t-s))\notag,
 \end{align}
In the following, we give the proofs for (\ref{A1_3}) and (\ref{A1_4}). The proof of (\ref{A1_5}) is similar to those for (\ref{A1_3}) and (\ref{A1_4}) and is not given here.

To prove (\ref{A1_3}) is true, we rewrite the right-hand side of (\ref{A1_3})
as double integrals
\small
\begin{align}
&\left( \int_0^t \sin(\omega (t-s)) \left(g(s)+  \frac{f(s)f_{\dot x}(s)}{2}\right) \,{\rm d}s \right)^2 \
 +   \left( \int_0^t \cos(\omega (t-s)) \left(g(s)+  \frac{f(s)f_{\dot x}(s)}{2}\right) \,{\rm d}s \right)^2\notag\\
 &=\int_0^t \int_0^t \sin(\omega (t-p))\sin(\omega (t-s))\left(g(p)+  \frac{f(p)f_{\dot x}(p)}{2}\right)\left(g(s)+  \frac{f(s)f_{\dot x}(s)}{2}\right) {\rm d}p\,{\rm d}s\notag\\
  &\quad +  \int_0^t \int_0^t \cos(\omega (t-p))\cos(\omega (t-s)) \left(g(p)+  \frac{f(p)f_{\dot x}(p)}{2}\right)\left(g(s)+  \frac{f(s)f_{\dot x}(s)}{2}\right) \,{\rm d}p\,{\rm d}s\notag\\
   &=\int_0^t \int_0^t \cos(\omega (s-p))\left(g(p)+  \frac{f(p)f_{\dot x}(p)}{2}\right)\left(g(s)+  \frac{f(s)f_{\dot x}(s)}{2}\right) \,{\rm d}p\,{\rm d}s\notag \\
    &=2\int_0^t \int_0^s \cos(\omega (s-p))\left(g(p)+  \frac{f(p)f_{\dot x}(p)}{2}\right)\left(g(s)+  \frac{f(s)f_{\dot x}(s)}{2}\right) \,{\rm d}p\,{\rm d}s\notag.
\end{align}
\normalsize
Similarly, to prove (\ref{A1_4}), we rewrite the right-hand side (\ref{A1_4}) as
\begin{align}\label{A1_8}
&\left( \int_0^t \sin(\omega (t-s)) f(s) \,{\rm d}B(s) \right)^2 + \left( \int_0^t \cos(\omega (t-s)) f(s) \,{\rm d}B(s) \right)^2 \nonumber\\
&=\int_0^t \int_0^t \cos(\omega (s-p)) f(s)f(p) \,{\rm d}B(p) \,{\rm d}B(s).
\end{align}
The integral domain for the right-hand side of (\ref{A1_8}) is a square given by  $A=\{(s,p) \big| s\in [0,t], p\in [0, t]\}$. Decompose the square into three parts: $A_1=\{(s,p)\big| 0\le s<p \le t\}$, $A_2=\{(s,p)\big| 0\le p<s\le t\}$, and $A_3=\{(s,s) \big| 0\le s \le t\}$, then the right-hand side
of (\ref{A1_8}) becomes
\begin{align}
 \int_0^t \int_0^t \cos(\omega (s-p)) f(s)f(p) \,{\rm d}B(p) \,{\rm d}B(s)
=\iint\limits_{A_1+A_2+A_3} \cos(\omega (s-p)) f(s) f(p) \,{\rm d}B(p) \,{\rm d}B(s).
\end{align}
Note that
\begin{align}\label{A1_9}
& \iint_{A_1} \cos(\omega (s-p)) f(s) f(p)  \,{\rm d}B(p) \,{\rm d}B(s) =
  \iint_{A_2} \cos(\omega (s-p)) f(s)f(p)  \,{\rm d}B(p) \,{\rm d}B(s) \nonumber\\
&= \int_0^t \int_0^s \cos(\omega (s-p)) f(s)f(p)  \,{\rm d}B(p) \,{\rm d}B(s) ,
\end{align}
and
\begin{align}\label{A1_10}
  \iint_{A_3} \cos(\omega (s-p)) f(s)f(p) \,{\rm d}B(p) \,{\rm d}B(s) =\int_0^t f^2(s) \,{\rm d}s .
 \end{align}
 It follows from (\ref{A1_9}) and (\ref{A1_10}) that (\ref{A1_4}) is true.

Add Eqs. (\ref{A1_3}) to  (\ref{A1_7}) together, we get $LHS=RHS$, and hence (\ref{energy_stratonovich_2}) is true.


\begin{thebibliography}{10}

\bibitem{DiPaolaFalsone1993}
M.~Di~Paola and G.~Falsone.
\newblock Ito and Stratonovich integrals for delta-correlated processes.
\newblock {\em Probabilistic engineering mechanics}, 8, 1993.

\bibitem{DiPaolaFalsone1993b}
M.~Di~Paola and G.~Falsone.
\newblock Stochastic dynamics of non-linear systems driven by non-normal
  delta-correlated processes.
\newblock {\em ASME Journal of applied mechanics}, 60:141--148, 1993.

\bibitem{Grigoriu1998}
M.~Grigoriu.
\newblock The Ito and Stratonovich integrals for stochastic differential
  equations with Poisson white noise.
\newblock {\em Probabilistic engineering mechanics}, 13:175--182, 1998.

\bibitem{Hu1994}
S.~L.~J. Hu.
\newblock Closure on discussion by Di Paola, m. and Falsone, g., on "response
  of dynamic systems excited by non-Gaussian pulse processes".
\newblock {\em ASCE Journal of engineering mechanics}, 120:2472--2474, 1994.

\bibitem{Ibrahim1985}
R.~A. Ibrahim.
\newblock {\em Parametric Random Vibration}.
\newblock Research Studies Press, 1985.

\bibitem{klebaner2005}
F.~C. Klebaner.
\newblock {\em Introduction to stochastic calculus with applications}.
\newblock 2nd Edition, Imperial College Press, 2005.

\bibitem{KloedenPlaten1992}
P.~Kloeden and E.~Platen.
\newblock {\em Numerical Solutions of Stochastic differential equations}.
\newblock Springer, 1992.

\bibitem{LinCai2004}
Y.~K. Lin and G.~Q. Cai.
\newblock {\em Probabilistic Structural Dynamics: Advanced Theory and
  Applications}.
\newblock Springer, 2005.

\bibitem{okse2003}
B.~K. Oksendal.
\newblock {\em Stochastic Differential Equations : an Introduction with
  Applications}.
\newblock Springer, 6th Edition, 2003.
%

\bibitem{CWSTo2000}
C.~W.~S. To.
\newblock {\em Nonlinear random vibration: Analytical techniques and
  applications}.
\newblock Swets and Zeitlinger Publishers, 2000.

\bibitem{CWSTo1988}
Cho W.~S. To.
\newblock On dynamic systems distributed by random parametric excitations.
\newblock {\em Journal of Sound and Vibration}, 123:387--390, 1988.

\bibitem{WongZakai1965}
E.~Wong and M.~Zakai.
\newblock On the relation between ordinary and stochastic differential
  equations.
\newblock {\em International Journal of Engineering Science}, 3:213--229, 1965.

\bibitem{YongLin1987}
Y.~Yong and Y.~K. Lin.
\newblock Exact stationary response solution for second order nonlinear systems
  under parametric and external white noise excitations.
\newblock {\em Journal of Applied Mechanics, Transactions of the American
  Society of Mechanical Engineers}, 54:414--418, 1987.

\end{thebibliography}
\end{document}